\documentclass[a4paper,12pt]{article}
\usepackage{amsmath}
\usepackage{amsthm}
\usepackage{amssymb}
\usepackage{amscd}

\title{Clifford Homomorphisms and Higher Spin Dirac Operators}
\author{Yasushi Homma\thanks{Department of Mathematics, 
  Waseda University, 
  Tokyo, Japan, \endgraf 
  {\it E-mail address}: homma@gm.math.waseda.ac.jp}}
\date{}

\theoremstyle{plain}
\newtheorem{thm}{Theorem}[section]
\newtheorem{prop}[thm]{Proposition}
\newtheorem{cor}[thm]{Corollary}
\newtheorem{lem}[thm]{Lemma}
\theoremstyle{definition}
\newtheorem{defini}{Definition}[section]

\numberwithin{equation}{section}

\theoremstyle{remark}
\newtheorem{rem}{Remark}[section]
\newtheorem{ex}{Example}[section]

\newcommand{\End}{\mathrm{End}}

\newcommand{\Ad}{\mathrm{Ad}}

\allowdisplaybreaks[3]
\begin{document}
\maketitle
%%%%%%%%%%%%%%%%%%%%%%%%%%%%%%
%        abstract            %
%%%%%%%%%%%%%%%%%%%%%%%%%%%%%%
\begin{abstract}
We present a generalization of the Clifford action for other representations spaces of $Spin(n)$, which is called the Clifford homomorphism. Their properties extend to the ones for the higher spin Dirac operators on spin manifolds. In particular, we have general Bochner identities for them, and an eigenvalue estimate of a Laplace type operator on any associated bundle.
\end{abstract}
%%%%%%%%%%%%%%%%%%%%%%%%%%%%%%%%
%%            0               %%
%%%%%%%%%%%%%%%%%%%%%%%%%%%%%%%%
\section{Introduction}\label{sec:0}
Let $M$ be a $n$-dimensional spin manifold and $\mathbf{Spin}(M)$ be its spin structure. The irreducible unitary representation $(\pi_{\rho},V_{\rho})$ of $Spin(n)$ induces a vector bundle $\mathbf{S}_{\rho}$ associated to $\mathbf{Spin}(M)$. We consider the canonical covariant derivative $\nabla$ mapping from $\Gamma(M,\mathbf{S}_{\rho})$ to $\Gamma(M,\mathbf{S}_{\rho}\otimes T^{\ast}(M))$. If we decompose the tensor bundle $\mathbf{S}_{\rho}\otimes T^{\ast}(M)$ into irreducible bundles, then we can construct first order differential operators associated to $\nabla$, 
\begin{equation}
D^{\rho}_{\lambda_k}:\Gamma (M,\mathbf{S}_{\rho})\xrightarrow{\nabla} \Gamma (M,\mathbf{S}_{\rho}\otimes T^{\ast}(M))\xrightarrow{\Pi^{\rho}_{\lambda_k}}\Gamma (M,\mathbf{S}_{\lambda_k}).  \label{eqn:0-1}
\end{equation}
Here $\Pi^{\rho}_{\lambda_k}$ is the orthogonal projection onto an irreducible bundle $\mathbf{S}_{\lambda_k}$ from $\mathbf{S}_{\rho}(M)\otimes T^{\ast}(M)$. These operators are called \textit{the higher spin Dirac operators}, \textit{the generalized gradient}, or \textit{the Stein-Weiss operators}, which are conformal invariant first order differential operators. Many basic and geometric operators on $M$ are constructed in this way; the Dirac operator, the twistor operator (see \cite{BFGK}, \cite{BF}, and \cite{LM}), the exterior derivative, its adjoint, the conformal killing operator (see \cite{GM} and \cite{K}), the Rarita-Schwinger operator (see \cite{Bu}) and so on. As in these examples, the higher spin Dirac operators and their properties are closely related with geometry and topology of $M$ (see the references given above and \cite{AHS}, \cite{Br}, \cite{Br2}, \cite{BS}, \cite{H3} etc.). 

Now, we consider the Dirac operator $D$ as a basic example. We denote the spinor representation of $Spin(n)$ by $(\pi_{\Delta},V_{\Delta})$. Then the Dirac operator $D$ is nothing else but the higher spin Dirac operator $D^{\Delta}_{\Delta}$, and the projection $\Pi^{\Delta}_{\Delta}$ is realized by using the Clifford algebra. The Clifford algebra gives us a lot of information about the Dirac operator. For instance, we can prove the ellipticity, the conformal invariance, and the Bochner identity of $D$ (see \cite{BFGK} and \cite{LM}). As this example, the projection $\Pi^{\rho}_{\lambda_k}$ in \eqref{eqn:0-1} is an essential tool to investigate the higher spin Dirac operators. Since the projection is defined fiberwise, we consider the tensor representation $(\pi_{\rho}\otimes \pi_{\Ad}, V_{\rho}\otimes \mathbf{R}^n)$ and the projection $\Pi^{\rho}_{\lambda_k}$ onto an irreducible component $V_{\lambda_k}$. For any $u$ in $\mathbf{R}^n$, we have a homomorphism $p^{\rho}_{\lambda_k}(u)$ from $V_{\rho}$ to $V_{\lambda_k}$,  
\begin{equation}
\mathbf{R}^n\times V_{\rho}\ni (u,\phi)\mapsto p^{\rho}_{\lambda_k}(u)\phi:=\Pi^{\rho}_{\lambda_k}(\phi \otimes u)\in V_{\lambda_k}.   
\end{equation}
We call this homomorphism \textit{the Clifford homomorphism}. On the spinor space, $p^{\Delta}_{\Delta}(u)$ is the usual Clifford action $u\cdot$. So the Clifford homomorphism is a generalization of the usual Clifford action.

The aim of this paper is to study the Clifford homomorphisms and the higher spin Dirac operators. Since we can extend the Clifford homomorphisms to bundle homomorphisms, some relations among the Clifford homomorphisms induce the ones among the higher spin Dirac operators. In particular, we have general Bochner identities. Furthermore, these identities allow us to give a lower bound of the first eigenvalue for a Laplace type operator on each associated bundle. 

In section \ref{sec:1}, we review the Clifford algebras, the spin groups, and their representations. In particular, we introduce the Casimir operator and the conformal weights. In section \ref{sec:2}, we decompose the representation space $V_{\rho}\otimes \mathbf{R}^n$ into irreducible components as a $Spin(n)$-module (or $\mathfrak{spin}(n)$-module) and define the Clifford homomorphisms. We investigate these homomorphisms and give an explicit formula of the projection $\Pi^{\rho}_{\lambda_k}$. As an example, we give explicit decompositions of $V_{\Delta}\otimes \mathbf{R}^n$ and $\Lambda^k(\mathbf{R}^n)\otimes \mathbf{R}^n$, where $(\pi_{\Delta},V_{\Delta})$ means the spinor representation and $(\pi_{\Lambda^k},\Lambda^k(\mathbf{R}^n))$ is the $k$-th exterior product representation. In section \ref{sec:3} and \ref{sec:4}, we study the higher spin Dirac operators and give some relations among them. Then we have the Bochner identities for the higher spin Dirac operators. In the last section, we give a lower bound of the first eigenvalue of a Laplace type operator constructed in section \ref{sec:4}. This is a generalization of the eigenvalue estimates known about the Dirac operator and the Laplace-Beltrami operator.
%%%%%%%%%%%%%%%%%%%%%%%%%%%%%%%%
%%           1                %%
%%%%%%%%%%%%%%%%%%%%%%%%%%%%%%%%
\section{Preliminaries: representations of $Spin(n)$}\label{sec:1}
In this section, we give a short review to the Clifford algebras, the spin groups, and their representations. Let $\mathbf{R}^n$ be the n-dimensional Euclidean space with orthonormal basis $\{e_k\}_{k=1}^{n}$. The Clifford algebra $Cl_n$ associated to $\mathbf{R}^n$ is an associative algebra with unit generated by $\{e_k\}$ under the relations $e_k e_l+e_l e_k=-2\delta_{kl}$. We denote the complexification of $Cl_n$ as $\mathbf{C}l_n$. 

We define the Lie algebra $\mathfrak{spin}(n)$ in $Cl_n$ by
%%%%%%%%%%%%% eqn 1-5 %%%%%%%
\begin{equation}
 \mathfrak{spin}(n):=\mathrm{span}_{\mathbf{R}} \{ [e_k, e_l] \:|\: 1\le k,l \le n\},
                                                            \label{eqn:1-2}
\end{equation}
%%%%%%
where $[e_k,e_l]:=e_ke_l-e_l e_k$ and the Lie bracket is $[a, b]:=ab-ba$ for $a$ and $b$ in $\mathfrak{spin}(n)$. The basis of $\mathfrak{spin}(n)$ is $\{[e_k,e_l]\}_{k<l}=\{2e_ke_l\}_{k<l}$. We put $\exp (a):=\sum\frac{a^n}{n!}$ for $a$ in $Cl_n$ and obtain the spin group $Spin(n)$, 
%%%%%%%%% eqn 1-6 %%%%%%%%%
\begin{equation}
Spin(n):=\exp \mathfrak{spin}(n) \subset Cl_n.         \label{eqn:1-3}
\end{equation}
If we would like to think of the spin group as a double covering group of $SO(n)$, we use the adjoint representation $\pi_{\Ad}$ of $Spin(n)$ on $\mathbf{R}^n$, 
%%%%%%% eqn  %%%%%%%%%%
\begin{equation}
Spin(n)\times \mathbf{R}^n \ni (g,x)\mapsto \pi_{\Ad}(g)x=gxg^{-1}\in \mathbf{R}^n.
                                \label{eqn:1-4}
\end{equation}
Then we have the isomorphisms $\mathfrak{so}(n)\simeq \mathfrak{spin}(n)$, where $[e_k,e_l]$ corresponds to the $n\times n$ matrix $-4E_{kl}+4E_{lk}$. The Killing form gives an inner product on $\mathfrak{spin}(n)$ as follows:
%%%%%%%%%%%%%
\begin{equation}
\langle [e_k,e_l],[e_i,e_j]\rangle =32\delta_{ik}\delta_{jl} \quad \textrm{for $k<l$ and $i<j$}.  \label{eqn:1-5}
\end{equation}
%%%%%%

The irreducible unitary representations of $Spin(n)$ or $\mathfrak{spin}(n)$ are parametrized by dominant weights $\rho=(\rho^1,\cdots,\rho^m) \in \mathbf{Z}^m \cup (\mathbf{Z}+1/2)^m$, satisfying that
%%%%%%%% weight %%%%%%%
\begin{gather}
\rho^1\ge \cdots \ge \rho^{m-1}\ge |\rho^m|, \quad \textrm{for $n=2m$}, \label{eqn:1-6} \\
\rho^1\ge \cdots \ge \rho^{m-1}\ge \rho^m\ge 0, \quad \textrm{for $n=2m+1$}.\label{eqn:1-7}
\end{gather}
%%%%%%
The dominant weight $\rho$ is the highest weight of the corresponding representation $(\pi_{\rho},V_{\rho})$. Here, we use the same notation for the representation of $Spin(n)$ and its infinitesimal representation of $\mathfrak{spin}(n)$. When writing dominant weights, we denote a string of $j$ $k$'s for $k$ in $\mathbf{Z}\cup (\mathbf{Z}+1/2)$ by $k_j$. For example, $((\frac{3}{2})_p,(\frac{1}{2})_{m-p})$ is the weight whose first $p$ components are $\frac{3}{2}$ and others are $\frac{1}{2}$. 

We shall introduce the Casimir operator and the conformal weights. Let $(\pi_{\rho},V_{\rho})$ be an irreducible representation with highest weight $\rho$. The Casimir operator on $V_{\rho}$ is defined by 
%%%%%%%%%% casimir %%%%%%%%%%%
\begin{equation}
C_{\pi_{\rho}}=\frac{1}{32}\sum_{i<j}\pi_{\rho}([e_i,e_j])\pi_{\rho}([e_i,e_j])=\frac{1}{64}\sum_{i,j}\pi_{\rho}([e_i,e_j])\pi_{\rho}([e_i,e_j]).
                \label{eqn:1-8}
\end{equation}
%%%%%%%
This operator commutes with the action of $Spin(n)$. So the Casimir operator is a constant $c(\rho)$ on $V_{\rho}$: 
%%%%%%% casimir scalar %%%%%
\begin{equation}
c(\rho):=-\frac{1}{2}(\|\delta+\rho\|^2-\|\delta\|^2). \label{eqn:1-9}
\end{equation}
%%%%
Here the inner products on the weight space is the standard one, that is, $\|\rho\|^2=\langle \rho,\rho\rangle:=\sum \rho^k\rho^k$, and $\delta$ is half the sum of the positive roots. 

We consider the tensor representation $(\pi_{\rho} \otimes \pi_{\Ad},V_{\rho}\otimes \mathbf{R}^n)\simeq (\pi_{\rho}\otimes \pi_{\Ad}, V_{\rho}\otimes \mathbf{C}^n)$ and its irreducible decomposition $\pi_{\rho}\otimes \pi_{\Ad}\simeq \pi_{\lambda_0}\oplus \pi_{\lambda_1} \oplus \cdots \oplus \pi_{\lambda_N}$. To identify the irreducible components, we need the Casimir operator or the operator $\widehat{C}$ given by
%%%%%%%%% def of C-hat %%%%%%%
\begin{equation}
\widehat{C}:=C_{\pi_{\rho}\otimes \pi_{\Ad}}
  -C_{\pi_{\rho}}\otimes 1-1\otimes C_{\pi_{\Ad}}. \label{eqn:1-12}
  \end{equation}
%%%%
It is easy to show that 
%%%%
\begin{equation}
\widehat{C}=\frac{1}{32}\sum_{i,j} \pi_{\rho}([e_i,e_j])\otimes \pi_{\Ad}([e_i,e_j]).
                              \label{eqn:1-13}
\end{equation}
%%%%
This operator $\widehat{C}$ is a constant on the irreducible component $V_{\lambda_k}$. This constant is called the conformal weight for $\lambda_k$ (see \cite{F}) and given by 
%%%%%%%% conformal weight %%%%%%%
\begin{equation}
\begin{split}
m(\lambda_k): &=c(\lambda_k)-c(\rho)-c(\Ad) \\
 &=\frac{1}{2}(n-\|\delta+\lambda_k\|^2+\|\delta+\rho\|^2-1). \label{eqn:1-14}
\end{split}
\end{equation}
%%%%%%
%%%%%%%%%%%%%%%%%%%%%%%%%%%%%%%%%%%%%%%%%%%%%%%%%%%%%%%%%%%%%%
%%%%%%%%%%%%%%%%%%%%%%%%%%%%%%%%
%%            2               %%
%%%%%%%%%%%%%%%%%%%%%%%%%%%%%%%%
\section{Clifford homomorphisms}\label{sec:2}
We shall discuss the Clifford homomorphisms. Before we define the Clifford homomorphisms, we recall the usual Clifford action on spinor spaces. We consider the complex spinor representation $(\Pi_{\Delta}, V_{\Delta})$ of $Spin(n)$. For $n=2m$, we remark that the spinor space $V_{\Delta}$ splits to $V_{\Delta^+}\oplus V_{\Delta^-}$ as a $Spin(n)$-module.  We adopt the following definition of the Clifford action. When $n$ is $2m+1$, the representation space $V_{\Delta}\otimes \mathbf{R}^n$ splits to $V_{T}\oplus V_{\Delta}$ as a $Spin(n)$-module, where the highest weights $\Delta$ and $T$ are $((1/2)_m)$ and $(3/2,(1/2)_{m-1})$ respectively. If we have a vector $\phi \otimes u$ in $V_{\Delta}\otimes \mathbf{R}^n$, then we can project the vector onto $V_{\Delta}$ along $V_T$ orthogonally. This projection gives the Clifford action of $\mathbf{R}^n$, $u\cdot \phi$. This action satisfy the relation $e_ie_j+e_je_i=-2\delta_{ij}$ and extends to the action of $\mathbf{C}l_n$. When $n$ is $2m$, the same discussion holds, so that we can define the Clifford action satisfying $\mathbf{R}^n(V_{\Delta^{\pm}})=V_{\Delta^{\mp}}$. We apply this definition to other representations. 

Let us consider an irreducible unitary representation $(\pi_{\rho},V_{\rho})$ and the tensor representation $\pi_{\rho}\otimes \pi_{\Ad}\simeq \pi_{\lambda_0}\oplus \cdots\oplus \pi_{\lambda_N}$. Since all weights of the adjoint representation $\pi_{\Ad}$ have multiplicity one, each irreducible constitute of $(\pi_{\rho}\otimes \pi_{\Ad}, V_{\rho}\otimes \mathbf{R}^n)$ has multiplicity one. In fact, Fegan shows the following fact in \cite{F}. 
%%%%%%%%%%%%%%%%%%%%%%%%%%%%%%%%%%%%%%%%%%%%%%%
%%%%%%%%%%%%%%%%  Lemma of Fegan     %%%%%%%%%%%%%%%%%%
%%%%%%%%%%%%%%%%%%%%%%%%%%%%%%%%%%%%%%%%%%%%%%%
\begin{lem}[\cite{F}]\label{lemma=F}
If $V_{\rho}\otimes \mathbf{R}^n=\sum V_{\lambda}$ is the irreducible decomposition as a $Spin(n)$-module (or a $\mathfrak{spin}(n)$-module), then the highest weight of irreducible components is given as follows:
\begin{enumerate}
\item When $n$ is $2m$, $\lambda$ is dominant weight and $\lambda=\rho\pm \mu_i$ $(1\le i \le m)$.
\item When $n$ is $2m+1$ and $\rho^m\ge\frac{1}{2}$, $\lambda$ is dominant weight, and $\lambda=\rho$ or $\lambda=\rho\pm \mu_i$ $(1\le i \le m)$. 
\item When $n$ is $2m+1$ and $\rho^m$ is $0$, $\lambda$ is dominant integral, and $\lambda=\rho+\mu_m$ or $\lambda=\rho\pm \mu_i$ $(1\le i \le m-1)$.
\end{enumerate}
Here $\rho=(\rho^1,\rho^2,\cdots,\rho^m)$ and $\mu_i=(0_{i-1},1,0_{m-i})$. The conformal weight of irreducible components is given as follows: 
%%%%%%%%%%%%
\begin{equation}
m(\rho+\mu_i)=i-1-\rho^i, \quad m(\rho-\mu_i)=n+\rho^i-i-1, \quad m(\rho)=\frac{1}{2}(n-1).                      \label{eqn:2-4}
\end{equation}
%%%%
\end{lem}
%%%%%%%%%%%%%%%%
We remark that the highest weight $\rho+\mu_1$ called \textit{top term} certainly occurs once in the decomposition. If we employ the lexicographical order for the weight space, then we arrange the highest weight $\{\lambda_k\}_{0\le k \le N}$ to satisfy $\rho+\mu_1=\lambda_0>\lambda_1>\cdots>\lambda_N$. Lemma \ref{lemma=F} implies that their conformal weights reverse as $-\rho^1=m(\lambda_0)<m(\lambda_1)<\cdots<m(\lambda_N)$ except the following case: when $n=2m$, $\rho^{m-1}\ge 1$, and $\rho^m=0$, we have $m(\rho+\mu_m)=m(\rho-\mu_m)$. This case is said to be \textit{the exceptional case}.

The inner products $\langle \cdot,\cdot \rangle$ on $V_{\rho}$ and $\mathbf{R}^n$ give the one on $V_{\rho}\otimes \mathbf{R}^n$ where irreducible components are orthogonal to each other. Then we have the orthogonal projection $\Pi^{\rho}_{\lambda_k}$ onto irreducible component $V_{\lambda_k}$ from $V_{\rho}\otimes \mathbf{R}^n$. This projection induces a homomorphism from $V_{\rho}$ to $V_{\lambda_k}$. 
%%%%%%%%%%%%%%%%%%%%%%%%%%%%%%%%%%%%%%%%%%%%%%%
%%%%%%%%%%%%%%%%% definition %%%%%%%%%%%%%%%%%%%%%
%%%%%%%%%%%%%%%%%%%%%%%%%%%%%%%%%%%%%%%%%%%%%%%
\begin{defini}\label{def}
Let $V_{\rho}\otimes \mathbf{R}^n=\sum V_{\lambda_k}$ be the irreducible decomposition as a $Spin(n)$-module. For each irreducible component $V_{\lambda_k}$, we have the bilinear mapping 
%%%%%%%%%% def of Clifford %%%%%%%
\begin{equation}
\mathbf{R}^n\times V_{\rho}\ni (u,\phi)\mapsto p^{\rho}_{\lambda_k}(u)\phi:=\Pi^{\rho}_{\lambda_k}(\phi \otimes u)\in V_{\lambda_k}.   \label{eqn:2-5}
\end{equation}
%%%%
We call the linear mapping $p^{\rho}_{\lambda_k}(u):V_{\rho}\to V_{\lambda_k}$ for all $u$ in $\mathbf{R}^n$, \textit{the Clifford homomorphism} from $V_{\rho}$ to $V_{\lambda_k}$.
\end{defini}
%%%%%%%%%%%%%%%%%%%%%%%%%%%%%%%%%%%
We denote by $(p^{\rho}_{\lambda_k}(u))^{\ast}$ the adjoint operator of $p^{\rho}_{\lambda_k}(u)$ with respect to the inner products on $V_{\rho}$ and $V_{\lambda_k}$. If we consider the tensor representation $V_{\lambda_k}\otimes \mathbf{R}^n =\sum V_{\nu}$, then we find the irreducible component $V_{\rho}$ and the Clifford homomorphism $p^{\lambda_k}_{\rho}(u)$ is $(p^{\rho}_{\lambda_k}(u))^{\ast}$ up to normalization.

From now on, we shall investigate the Clifford homomorphisms. Since the Clifford homomorphisms is defined by a projection mapping, we can easily notice the following fact.
%%%%%%%%%%%%%%%%%%%%%%%%%%%%%%%%%%%%%%%%%%%%%%%
%%%%%%%%%%%% proposition (relation 1) %%%%%%%%%%%%%%%%%%%
%%%%%%%%%%%%%%%%%%%%%%%%%%%%%%%%%%%%%%%%%%%%%%%
\begin{prop}\label{relation 1}
For $u,v$ in $\mathbf{R}^n$, the Clifford homomorphisms $(p^{\rho}_{\lambda_k}(u))^{\ast}p^{\rho}_{\lambda_k}(v)$ satisfy that  
%%%%
\begin{equation}
\sum_k (p^{\rho}_{\lambda_k}(u))^{\ast}p^{\rho}_{\lambda_k}(v)=\langle u,v\rangle \quad \textrm{on $V_{\rho}$}.          \label{eqn:2-6}
\end{equation}
%%%%
\end{prop}
%%%%%%%%%%%%%
%%%%%%%%%%%%%%%% begin Proof %%%%%%%%%%%%%%%%%%%%%%
\begin{proof}
For $\phi\otimes v$ and $\psi\otimes v$ in $V_{\rho}\otimes \mathbf{R}^n$, we have 
$$
\langle \phi \otimes v, \psi \otimes u \rangle=\langle \phi,\psi \rangle \langle u,v\rangle=\langle \langle u,v \rangle \phi,\psi\rangle.
$$
On the other hand, we have 
\begin{equation}
\langle \phi \otimes v, \psi \otimes u\rangle = \sum_k \langle p^{\rho}_{\lambda_k}(v)\phi, p^{\rho}_{\lambda_k}(u)\psi \rangle = \langle \sum_k ( p^{\rho}_{\lambda_k} (u) )^{\ast} p^{\rho}_{\lambda_k} (v) \phi, \psi \rangle. \nonumber 
\end{equation}
Since these equations hold for all $\psi$ in $V_{\rho}$, we have proved the proposition.
\end{proof}
%%%%%%%%%%%%%%% end proof %%%%%%%%%%%%%%%%%
Next, we show an important lemma to give relations among the Clifford homomorphisms.  
%%%%%%%%%%%%%%%%%%%%%%%%%%%%%%%%%%%%%%%%%%%%%%%
%%%%%%%%%%%% lemma %%%%%%%%%%%%%%%%%%%%%%%%%%%%%%%%%%%
%%%%%%%%%%%%%%%%%%%%%%%%%%%%%%%%%%%%%%%%%%%%%%%
\begin{lem}\label{lemma 1}
The Clifford homomorphism $p^{\rho}_{\lambda_k}(u)$ from $V_{\rho}$ to $V_{\lambda_k}$ satisfies the relation 
%%%%%%
\begin{equation}
-\frac{1}{4} \sum_i p^{\rho}_{\lambda_k}(e_i)\pi_{\rho}([e_i,u])=m(\lambda_k)p^{\rho}_{\lambda_k}(u),   \label{eqn:2-7}
\end{equation}
%%%%
where $\{e_i\}_{i=1}^n$ is any orthonormal basis of $\mathbf{R}^n$.
\end{lem}
%%%%%%%%%%%%%%%
%%%%%%%%%%%%%% begin proof %%%%%%%%%%%%%%%%%%%%%%%%
\begin{proof}
To prove the lemma, we use the operator $\widehat{C}$ given in section \ref{sec:1}. For $\phi\otimes u$ in $V_{\rho}\otimes \mathbf{R}^n$, we have 
%%%%
\begin{equation}
\begin{split}
\widehat{C}(\phi\otimes u)&=\frac{1}{32}\sum_{s,t,i}\pi_{\rho}([e_s,e_t])\phi\otimes [[e_s,e_t], e_i ] \langle e_i,u \rangle \\ 
&=\frac{1}{8}\sum_{s,t} \pi_{\rho}([e_s,e_t])\phi\otimes (\delta_{si}e_t-\delta_{ti}e_s) \langle e_i,u\rangle \\ 
&=\frac{1}{4}\sum_l \pi_{\rho}([u,e_t])\phi \otimes e_t \\
&=\frac{1}{4}\sum_{k,t} p^{\rho}_{\lambda_k}(e_t)\pi_{\rho}([u,e_t])\phi. 
\end{split} \nonumber
\end{equation}
%%%
The operator $\widehat{C}$ is the constant $m(\lambda_k)$ on $V_{\lambda_k}$, so that $\widehat{C}(\phi \otimes u)=\sum_k m(\lambda_k)p^{\rho}_{\lambda_k}(u)\phi$. Hence, we conclude that 
%%%
\begin{equation}
m(\lambda_k)p^{\rho}_{\lambda_k}(u)=\frac{1}{4}\sum_{t} p^{\rho}_{\lambda_k}(e_t)\pi_{\rho}([u,e_t]).   \nonumber
\end{equation}
\end{proof}
%%%%%%%%%%%%%%% end Proof %%%%%%%%%%%%%%%%%%%%%%%%%%%%%%%%%%%%
This lemma leads a relation among the Clifford homomorphisms.
%%%%%%%%%%%%%%%%%%%%%%%%%%%%%%%%%%%%%%%%%%%%%%%%%%%%%%%%%%%%%%%%%%
%%%%%%%%%%%%%%%%%%%%% proposition (relation 2) %%%%%%%%%%%%%%%%%%%%
%%%%%%%%%%%%%%%%%%%%%%%%%%%%%%%%%%%%%%%%%%%%%%%%%%%%%%%%%%%%%%%%%%
\begin{prop}\label{relation 2}
We sum up the Clifford homomorphisms $(p^{\rho}_{\lambda_k}(u))^{\ast}p^{\rho}_{\lambda_k}(v)$ for all $\lambda_k$ with its conformal weight. Then we have
%%%
\begin{equation}
\sum_k m(\lambda_k)(p^{\rho}_{\lambda_k}(u))^{\ast}p^{\rho}_{\lambda_k}(v)=-\frac{1}{4}\pi_{\rho}([u,v]). \label{eqn:2-8}
\end{equation}
%%%
\end{prop}
%%%%%%%%%
%%%%%%%%%%%%%%%%%%%% begin proof %%%%%%%%%%%%%%%%%%%%
\begin{proof}
We substitute \eqref{eqn:2-7} into $\sum_{k}m(\lambda_k)(p^{\rho}_{\lambda_k}(u))^{\ast}p^{\rho}_{\lambda}(v)$ and use the relation \eqref{eqn:2-6}. 
\end{proof}
%%%%%%%%%%%%%% end proof %%%%%%%%%%%%%%%%%%%%
As this proposition, we give many relations among the Clifford homomorphisms.
%%%%%%%%%%%%%%%%%%%%%%%%%%%%%%%%%%%%%%%%%%%%%%%%%%%%%%%%%%%%%%%%%%
%%%%%%%%%%%%%%%%%%%%%%% theorem (relation for n) %%%%%%%%%%%%%%%
%%%%%%%%%%%%%%%%%%%%%%%%%%%%%%%%%%%%%%%%%%%%%%%%%%%%%%%%%%%%%%%%%%
\begin{thm}\label{relation n}
For an irreducible representation $(\pi_{\rho}, V_{\rho})$ and any non-negative integer $q$, 
we define the bilinear mapping $r^q_{\rho}$ from $\mathbf{R}^n\times \mathbf{R}^n$ to $\End (V_{\rho})$ as follows: 
\begin{multline}
r^{q}_{\rho}:\mathbf{R}^n\times \mathbf{R}^n \ni (u,v) \mapsto \\
    \left( -\frac{1}{4} \right)^{q} \sum_{l_1,\cdots, l_{q-1}} \pi_{\rho}([u, e_{l_1}])\pi_{\rho}([e_{l_1},e_{l_2}]) \cdots \pi_{\rho}([e_{l_{q-1}}, v])\in \mathrm{End}(V_{\rho}), \label{eqn:2-9}
 \end{multline}
and $r^0_{\rho}(u,v):=\langle u,v \rangle$. Then we have relations among the Clifford homomorphisms for any $q$: 
\begin{equation}
\sum_k m(\lambda_k)^{q}(p^{\rho}_{\lambda_k}(u))^{\ast}p^{\rho}_{\lambda_k}(v)=r^q_{\rho}(u,v).   \label{eqn:2-10}
\end{equation}
\end{thm}
%%%%%%%%%%%%%%%%%
As a corollary of this theorem, we give a description of general Casimir operators. Let $U(\mathfrak{spin}(n))$ be the enveloping algebra of $\mathfrak{spin}(n)$. Then the general Casimir operator is said to be an element of the center of $U(\mathfrak{spin}(n))$ or its image by $\pi_{\rho}$ for the representation $(\pi_{\rho},V_{\rho})$. In \cite{Z}, we know generators of the general Casimir operators and its action on the irreducible representations. The generators are given by 
\begin{equation}
\left(-\frac{1}{4}\right)^{q} \sum_{l_1, \cdots ,l_{q}}\pi_{\rho}([e_{l_1}, e_{l_2}])\pi_{\rho}([e_{l_2},e_{l_3}])\cdots \pi_{\rho}([e_{l_{q}}, e_{l_1}]).
\label{eqn:2-11}
\end{equation}
%%%%%%%%%%%%%%%%%%%%%%%%%%%%%%%%%%%%%%%%%%%%%%%
%%%%%%%%%%%%%%%%% Corollary Casimir %%%%%%%%%%%%%%%
%%%%%%%%%%%%%%%%%%%%%%%%%%%%%%%%%%%%%%%%%%%%%%%
\begin{cor}\label{casimir}
The Clifford homomorphisms $\sum_{k,i} m(\lambda_k)^{q}(p^{\rho}_{\lambda_k}(e_i))^{\ast}p^{\rho}_{\lambda_k}(e_i)$ for $q\ge 0$ generate the general Casimir operators for $(\pi_{\rho},V_{\rho})$ except the Pfaffian type Casimir operator in the remark below.
\end{cor}
%%%%%%%%%%%%%%%%%%%%%%%%%%%%%%%%%%%%%%%%%%%%%%%
%%%%%%%%%%%%%%%%%%%%%%   remark %%%%%%%%%%%%%%%
%%%%%%%%%%%%%%%%%%%%%%%%%%%%%%%%%%%%%%%%%%%%%%%
\begin{rem}\label{remark:2-1}
For $n=2m$, we need the Pfaffian type Casimir operator to generate the center of $U(\mathfrak{spin}(2m))$, that is,
\begin{equation}
\mathrm{Pf}_{\rho}:=\sum_{\sigma \in S_{2m}} \mathrm{sig}(\sigma) \pi_{\rho}([e_{\sigma (1)}, e_{\sigma (2)}])\pi_{\rho}([e_{\sigma (3) },e_{\sigma (4)}]) \cdots \pi_{\rho}([e_{\sigma (2m-1)}, e_{\sigma (2m)}]). \label{eqn:2-12}
\end{equation}
We can show a relation between this operator and the Clifford homomorphisms, 
\begin{multline}
\sum_k p(\lambda_k)(p^{\rho}_{\lambda_k}(e_j))^{\ast}p^{\rho}_{\lambda_k}(e_i)=\delta_{ij}p(\rho)+ \\ 
 8m(1-\delta_{ij})\mathrm{sgn}
\left( \begin{smallmatrix}
1& 2 &\cdots &2m \\
i& j &\cdots &2m
\end{smallmatrix} \right)
\sum_{\sigma \in \widetilde{S}_{2m}} \mathrm{sgn}(\sigma) \pi_{\rho}([e_{\sigma (1)}, e_{\sigma (2)}])\cdots \pi_{\rho}([e_{\sigma (2m-1)}, e_{\sigma (2m)}]),  \label{eqn:2-13}
\end{multline}
where $\widetilde{S}_{2m}$ is the permutation of $\{ 1,\cdots ,2m \}\setminus \{ i,j \}$ and $p(\rho)$ is a scalar action of $\mathrm{Pf}_{\rho}$ on $V_{\rho}$ as 
\begin{equation}
p(\rho)=(\sqrt{-1})^m 8^m m!(\rho^1+m-1)(\rho^2+m-2)\cdots 
(\rho^{m-1}+1)\rho^m.  \label{eqn:2-14}
\end{equation}
If we define the bilinear mapping $\mathrm{pf}_{\rho}(\cdot,\cdot)$ from $\mathbf{R}^n\times \mathbf{R}^n$ to $\End (V_{\rho})$ by the right hand side of \eqref{eqn:2-13}, then the trace of $\mathrm{pf}_{\rho}(\cdot,\cdot)$ is $n \mathrm{Pf}_{\rho}$.
\end{rem}
%%%%%%%%%%%%%%%%%%%%  remark end   %%%%%%%%%%%%%%%%%%%%

Now, we know that the usual Clifford action satisfies that
\begin{equation}
(g u g^{-1})\cdot \phi=\pi_{\Delta}(g)(u\cdot(\pi_{\Delta}(g^{-1})\phi)), \label{eqn:2-16}
\end{equation}
where $\phi$ is in $V_{\Delta}$ and $g$ in $Spin(n)$ (see \cite{LM}). We generalize this relation. 
%%%%%%%%%%%%%%%%%%%%%%%%%%%%%%%%%%%%%%%%%%%%%%%
%%%%%%%%%%%%%% Proposition (relation 3) %%%%%%%%%%%%%%
%%%%%%%%%%%%%%%%%%%%%%%%%%%%%%%%%%%%%%%%%%%%%%%
\begin{prop}\label{relation 3}
Let $u$ be in $\mathbf{R}^n$ and $g$ be in $Spin(n)$. The Clifford homomorphism $p^{\rho}_{\lambda_k}(u)$ satisfy that 
\begin{equation}
p^{\rho}_{\lambda_k}(g u g^{-1})=\pi_{\lambda_k}(g)p^{\rho}_{\lambda_k}
(u)\pi_{\rho}(g^{-1}).         \label{eqn:2-17}
\end{equation}
\end{prop}
%%%%%%%%%%%%%
%%%%%%%%%%%%%%%%%%%%% begin proof %%%%%%%%%%%%%
\begin{proof}
For $\phi$ in $V_{\rho}$, we have 
\begin{equation}
\begin{split}
p^{\rho}_{\lambda_k}(g u g^{-1})\phi &=\Pi^{\rho}_{\lambda_k}(\pi_{\rho}(g)\pi_{\rho}(g^{-1})\phi\otimes \pi_{\Ad}(g)u) \\
&=\Pi^{\rho}_{\lambda_k}(\pi_{\rho\otimes\Ad}(g)(\pi_{\rho}(g^{-1})\phi\otimes u)) \\
&=\pi_{\lambda_k}(g)\Pi^{\rho}_{\lambda_k}(\pi_{\rho}(g^{-1})\phi\otimes u) \\
&=\pi_{\lambda_k}(g)p^{\rho}_{\lambda_k}
(u)\pi_{\rho}(g^{-1})\phi.
\end{split} \nonumber
\end{equation}
\end{proof}
%%%%%%%%%%%% end proof %%%%%%%%%%%%%%%%%%
As a corollary of Proposition \ref{relation 3}, we have a formula of the projection $\Pi^{\rho}_{\lambda_k}$ by using the Clifford homomorphisms.
%%%%%%%%%%%%%%%%%%%%%%%%%%%%%%%%%%%%%%%%%%%%%%%
%%%%%%%%%%%%%% Corollary %%%%%%%%%%%%%
%%%%%%%%%%%%%%%%%%%%%%%%%%%%%%%%%%%%%%%%%%%%%%%
\begin{cor}\label{projection}
The orthogonal projection $\Pi_{\lambda_k}^{\rho}:V_{\rho}\otimes \mathbf{R}^n\to V_{\lambda_k}$ is realized as follows:
\begin{equation}
\Pi_{\lambda_k}^{\rho}(\phi \otimes u)=\sum_i (p_{\lambda_k}^{\rho}(e_i))^{\ast}p^{\rho}_{\lambda_k}(u)(\phi)\otimes e_i.   \label{eqn:2-18}
\end{equation}
\end{cor}
%%%%%
%%%%%%%%%%% begin proof %%%%%%%%%%
\begin{proof}
Let $V_{\lambda_k}$ be an irreducible component of $V_{\rho}\otimes \mathbf{R}^n$. We consider the following embedding from $V_{\lambda_k}$ to $V_{\rho}\otimes \mathbf{R}^n$, which dose not always preserve their inner products:
\begin{equation}
i_{\lambda_k}:V_{\lambda_k}\ni \psi \mapsto \sum_i (p_{\lambda_k}^{\rho}(e_i))^{\ast}(\psi)\otimes e_i 
    \in V_{\rho}\otimes \mathbf{R}^n. \nonumber
\end{equation}
The map $i_{\lambda_k}$ is independent of the orthonormal basis of $\mathbf{R}^n$ which we chose, and commutes with the action of $Spin(n)$ . Since the irreducible component of $V_{\rho}\otimes \mathbf{R}^n$ has multiplicity one and $i_{\lambda_k}$ is not zero, we prove that the map $i_{\lambda_k}$ is a well-defined embedding. 

For $\phi\otimes u$ in $V_{\rho}\otimes \mathbf{R}^n$, we have 
\begin{equation}
\begin{split}
\phi \otimes u &= \sum_i \langle u,e_i\rangle \phi\otimes e_i \\
&=\sum_{k,i} (p_{\lambda_k}^{\rho}(e_i))^{\ast}p^{\rho}_{\lambda_k}(u)(\phi)\otimes e_i \quad \textrm{( from Proposition \ref{relation 1} )} \\
&=\sum_k i_{\lambda_k}(p^{\rho}_{\lambda_k}(u)(\phi)).
\end{split} \nonumber
\end{equation}
This completes the proof of Corollary \ref{projection}.
\end{proof}
%%%%%%%%%%%%% end proof %%%%%%%%%%%%%%%%%%%%%%%
This corollary and Theorem \ref{relation n} allow us to have an explicit formula of the projection $\Pi^{\rho}_{\lambda_k}$. We define the $(N+1) \times (N+1)$ matrix $\mathbf{M}=(m_{ij})_{0 \le i,j \le N}$ by $m_{ij}:=m(\lambda_j)^i$, which is a Vandermonde matrix. Since the conformal weights are different from each other except \textit{the exceptional case}, the inverse matrix of $\mathbf{M}$ exists. From \eqref{eqn:2-10} and \eqref{eqn:2-18}, we have a formula of $\Pi_{\lambda_k}^{\rho}$, 
\begin{equation}
\Pi_{\lambda_k}^{\rho}(\phi \otimes u)=\sum_{i,q} n_{k q}r^q_{\rho}(e_i,u)\phi\otimes e_i.   \label{eqn:2-19}
\end{equation}
Here, $n_{ij}$ is the $(i,j)$-component of $\mathbf{M}^{-1}$. If we set $S_j(x_0,\cdots,\hat{x_i},\cdots,x_N)$ as the $j$-th fundamental symmetric polynomial of $\{x_0,x_1,\cdots,x_N\}\setminus \{x_i\}$, then 
\begin{equation}
n_{ij}=(-1)^{N-j}\frac{S_{N-j}(m(\lambda_0),\cdots, \widehat{m(\lambda_i)},\cdots, m(\lambda_N))}{\prod_{k\neq i}(m(\lambda_i)-m(\lambda_k))}.  
       \label{eqn:2-20}
\end{equation}

In the exceptional case that $n$ is even, $\rho^m\ge 1$, and $\rho^m=0$, the conformal weight of $\lambda_+:=\rho+\mu_m$ coincides with the one of $\lambda_-:=\rho-\mu_m$. So we use $\mathrm{pf}_{\rho}$ given in \eqref{eqn:2-13} to obtain a formula of $\Pi_{\lambda_{\pm}}^{\rho}$. It is easy to show that 
\begin{equation}
p(\lambda_+)(p^{\rho}_{\lambda_+}(e_j))^{\ast}
      p^{\rho}_{\lambda_+}(e_i)+p(\lambda_-)(p^{\rho}_{\lambda_-}(e_j))^{\ast}
      p^{\rho}_{\lambda_-}(e_i)=\mathrm{pf}_{\rho}(e_j,e_i). \label{eqn:2-21}
\end{equation}
It follows that we realize the projection of $\Pi_{\lambda_{\pm}}^{\rho}$ by using $r_{\rho}^q$ and $\mathrm{pf}_{\rho}$.

We shall state some examples of the Clifford homomorphisms.
%%%%%%%%%%%%%%%%%%%%%%%%%%%%%%%%%%%%%%%%%%%%%%%%%%%%%%%%%%%
%%%%%%%%%%      example 1    %%%%%%%%%%%%%%%%%%%%%%%%%%%%%%
%%%%%%%%%%%%%%%%%%%%%%%%%%%%%%%%%%%%%%%%%%%%%%%%%%%%%%%%%%%
\begin{ex}[Spinor]\label{ex:spinor 1}
We discuss only the case of $n=2m+1$ and leave the case of $n=2m$ to the reader. Let $V_{\Delta}$ be spinor space which is an irreducible representation space with highest weight $((1/2)_{m})$. It follows from Lemma \ref{lemma=F} that the vector space $V_{\Delta}\otimes \mathbf{R}^n$ is isomorphic to $V_{\Delta}\oplus V_{T}$. Here, $V_T$ is the representation space with highest weight $(3/2,(1/2)_{m-1})$. From Theorem \ref{relation n}, we have 
\begin{gather}
(p^{\Delta}_{T}(u))^{\ast}p^{\Delta}_{T}(v)+(p^{\Delta}_{\Delta}(u))^{\ast}p^{\Delta}_{\Delta}(v)=\langle  u,v\rangle , \label{eqn:2-e-1} \\
-\frac{1}{2}(p^{\Delta}_{T}(u))^{\ast}p^{\Delta}_{T}(v)+\frac{n-1}{2}(p^{\Delta}_{\Delta}(u))^{\ast}p^{\Delta}_{\Delta}(v)=-\frac{1}{4}\pi_{\Delta}([u,v]), 
           \label{eqn:2-e-2}
\end{gather}
for $u$ and $v$ in $\mathbf{R}^n$. It follows that we give the usual Clifford relation 
\begin{equation}
(p^{\Delta}_{\Delta}(u))^{\ast}p^{\Delta}_{\Delta}(v)+(p^{\Delta}_{\Delta}(v))^{\ast}p^{\Delta}_{\Delta}(u)=\frac{2}{n}\langle u,v\rangle. \label{eqn:2-e-4}
\end{equation}
Since $p^{\Delta}_{\Delta}$ is a constant multiple of the Clifford action $u\cdot $ satisfying $\langle u\cdot \phi,\psi\rangle =-\langle \phi,u\cdot\psi\rangle$, we have $p^{\Delta}_{\Delta}(u)=\frac{1}{\sqrt{n}}u\cdot$ and $(p^{\Delta}_{\Delta}(u))^{\ast}=-\frac{1}{\sqrt{n}}u\cdot$. Furthermore the representation of $\mathfrak{spin}(n)$ on spinor space is realized as $\pi_{\Delta}([e_i,e_j])=[e_i\cdot,e_j\cdot]$. The orthogonal projections $\Pi_{\Delta}^{\Delta}$ and $\Pi_{T}^{\Delta}$ are 
\begin{gather}
\begin{split}
\Pi_{\Delta}^{\Delta}(\phi\otimes u) &=\sum \{\frac{1}{n}\langle e_i,u\rangle-\frac{1}{2n}\pi_{\Delta}([e_i,u])\}( \phi ) \otimes e_i \\
&= -\frac{1}{n}\sum e_i\cdot u\cdot \phi \otimes e_i,    
\end{split}\label{eqn:2-e-5} \\
\begin{split}
 \Pi_{T}^{\Delta}(\phi\otimes u)&=\sum \{\frac{n-1}{n}\langle e_i,u\rangle+\frac{1}{2n}\pi_{\Delta}([e_i,u])\}( \phi ) \otimes e_i \\
 &=\phi \otimes u+\frac{1}{n} \sum e_i\cdot u\cdot \phi \otimes e_i.
\end{split} \label{eqn:2-e-6}
\end{gather}
\end{ex}
%%%%%%%%%%%%%%%%%%%%%%%%%%%%%%%%%%%%%%%%%%%%%%%%%%%%%%%%%%%
%%%%%%%%%%      example 2    %%%%%%%%%%%%%%%%%%%%%%%%%%%%%%
%%%%%%%%%%%%%%%%%%%%%%%%%%%%%%%%%%%%%%%%%%%%%%%%%%%%%%%%%%%
\begin{ex}[Exterior algebra]\label{ex:diff 1}
We shall discuss the exterior tensor product representation $(\pi_{\Lambda^k},\Lambda^k(\mathbf{R}^n))$ with highest weight $\rho=(1_k,0_{m-k})$ for $1\le k \le m$. There are three irreducible representations $\lambda_0=(2,1_{k-1},0_{m-k})$, $\lambda_1=(1_{k+1},0_{m-k-1})$, and $\lambda_2=(1_{k-1},0_{m-k+1})$ in $\Lambda^k(\mathbf{R}^n)\otimes \mathbf{R}^n$. We realize the Clifford homomorphisms as follows:
\begin{align}
(p^{\rho}_{\lambda_0}(e_j))^{\ast}p^{\rho}_{\lambda_0}(e_i)&=\frac{1}{(k+1)(n-k+1)}(k(n-k)\delta_{ij}-n r^1_{\rho}(e_j,e_i)+r^2_{\rho}(e_j,e_i)), 
        \label{eqn:2-e-7}  \\
(p^{\rho}_{\lambda_1}(e_j))^{\ast}p^{\rho}_{\lambda_1}(e_i)&=\frac{1}{(k+1)(n-2k)}((n-k)\delta_{ij}+(n-k-1)r^1_{\rho}(e_j,e_i)-r^2_{\rho}(e_j,e_i)),
       \label{eqn:2-e-8} \\
(p^{\rho}_{\lambda_2}(e_j))^{\ast}p^{\rho}_{\lambda_2}(e_i)&=\frac{1}{(n-2k)(n-k+1)}(-k\delta_{ij}-(k-1)r^1_{\rho}(e_j,e_i)+r^2_{\rho}(e_j,e_i)).  
     \label{eqn:2-e-9}
\end{align}
In general, it holds that $r^2_{\rho}(u,v)-r^2_{\rho}(v,u)=(n-2)r^1_{\rho}(u,v)$. Then we have 
\begin{equation}
(k+1)(p^{\rho}_{\lambda_1}(e_j))^{\ast}p^{\rho}_{\lambda_1}(e_i)+(n-k+1)(p^{\rho}_{\lambda_2}(e_i))^{\ast}p^{\rho}_{\lambda_2}(e_j)  =\delta_{ij},\label{eqn:2-e-10} 
\end{equation}
and 
\begin{equation}
\begin{split}
 &(k+1)\{(p^{\rho}_{\lambda_1}(e_j))^{\ast}p^{\rho}_{\lambda_1}(e_i)-(p^{\rho}_{\lambda_1}(e_i))^{\ast}p^{\rho}_{\lambda_1}(e_j)\}  \\
=&(n-k+1)\{(p^{\rho}_{\lambda_2}(e_j))^{\ast}p^{\rho}_{\lambda_2}(e_i)-(p^{\rho}_{\lambda_2}(e_i))^{\ast}p^{\rho}_{\lambda_2}(e_j)\} \\
=&r^1_{\rho}(e_j,e_i)=-\frac{1}{4}\pi_{\rho}([e_j,e_i]).
\end{split}
     \label{eqn:2-e-10-2}
\end{equation}
These equations correspond to the following relations, respectively:
\begin{gather}
i(e_j)e_{i\wedge}+e_{j\wedge}i(e_i)=\delta_{ij},  \label{eqn:2-e-11} \\
i(e_i)e_{j\wedge}-i(e_i)e_{j\wedge}=e_{j\wedge}i(e_i)-e_{i\wedge}i(e_j)=-\frac{1}{4}\pi_{\rho}([e_j,e_i]),   \label{eqn:2-e-11-2}
\end{gather}
where $i(u)$ is the interior product of $u$. So we put 
\begin{equation}
p^{\rho}_{\lambda_1}(e_i):=\frac{1}{\sqrt{k+1}}e_{i\wedge}, \quad 
    p^{\rho}_{\lambda_2}(e_i):=\frac{1}{\sqrt{n-k+1}}i(e_i). \label{eqn:2-e-12}
\end{equation}
Then the projections are realized as 
\begin{gather}
\Pi^{\rho}_{\lambda_1}(\phi\otimes u)= \frac{1}{k+1}\sum i(e_i)u_{\wedge}\phi\otimes e_i,   \label{eqn:2-e-13} \\
\Pi^{\rho}_{\lambda_2}(\phi\otimes u)= \frac{1}{n-k+1}\sum e_{i\wedge}i(u)\phi\otimes e_i, \label{eqn:2-e-14} \\
\Pi^{\rho}_{\lambda_0}(\phi\otimes u)= \phi\otimes u-\Pi^{\rho}_{1}(\phi\otimes u)-\Pi^{\rho}_{2}(\phi\otimes u).  \label{eqn:2-e-15}
\end{gather} 
If we consider the (exceptional) case that $n$ is $2m$ and the highest weight $\rho$ is $(1_{m-1},0)$, then we have four irreducible representations, $\lambda_0=(2,1_{m-1},0)$, $(\lambda_1)_+=(1_m)$, $(\lambda_1)_-=(1_{m-1},-1)$, and $\lambda_2=(1_{m-2},0_2)$ in $\Lambda^{m-1}(\mathbf{R}^{2m})\otimes \mathbf{R}^{2m}$. We denote the Hodge star operator by the asterisk $\ast$ and have 
\begin{gather}
\Pi^{\rho}_{(\lambda_1)_+}(\phi\otimes u)= \frac{1}{m}\sum i(e_i)\frac{1}{2}(1+\ast) u_{\wedge}\phi\otimes e_i,   \label{eqn:2-e-16}\\
\Pi^{\rho}_{(\lambda_1)_-}(\phi\otimes u)= \frac{1}{m}\sum i(e_i)\frac{1}{2}(1-\ast)u_{\wedge}\phi\otimes e_i.  \label{eqn:2-e-17}
\end{gather}
\end{ex}
%%%%%%%%%%%%%%%%%%%%%%%%%%%%%%%
%%            3               %%
%%%%%%%%%%%%%%%%%%%%%%%%%%%%%%%%
\section{Higher spin Dirac operators}\label{sec:3}
In this section we shall discuss the higher spin Dirac operators. Let $M$ be a $n$-dimensional spin manifold and $\mathbf{Spin}(M)$  be its spin structure (If $M$ dose not have a spin structure, we should consider only the representations of $SO(n)$). For an irreducible unitary representation $(\pi_{\rho},V_{\rho})$ of the structure group $Spin(n)$, we have a vector bundle $\mathbf{S}_{\rho}$ associated to $\mathbf{Spin}(M)$, 
%%%%
\begin{equation}
\mathbf{S}_{\rho}:=\mathbf{Spin}(M)\times_{\rho}V_{\rho}.
                              \label{eqn:3-1}
\end{equation}
%%%
We introduce a covariant derivative on $\Gamma(M,\mathbf{S}_{\rho})$ associated to the spin connection or the Levi-Civita connection, 
%%%%
\begin{equation}
\nabla: \Gamma(M,\mathbf{S}_{\rho})\to \Gamma(M,\mathbf{S}_{\rho}\otimes T^{\ast}(M)).                            \label{eqn:3-2}
\end{equation}
%%%%
If we decompose $\mathbf{S}_{\rho}\otimes T^{\ast}(M)$ into the direct sum of irreducible components, $\sum \mathbf{S}_{\lambda_k}$, then we construct a first order differential operator from $\Gamma(M,\mathbf{S}_{\rho})$ to $\Gamma(M,\mathbf{S}_{\lambda_k})$ called the higher spin Dirac operator,
%%%%%%%%%%%%%%%
\begin{equation}
D^{\rho}_{\lambda_k}:\Gamma (M,\mathbf{S}_{\rho})\xrightarrow{\nabla} \Gamma (M,\mathbf{S}_{\rho}\otimes T^{\ast}(M))\xrightarrow{\Pi^{\rho}_{\lambda_k}}\Gamma (M,\mathbf{S}_{\lambda_k}).               \label{eqn:3-3}
\end{equation}
%%%%%%%%%%%
From Proposition \ref{relation 3}, we show that the bundle homomorphism $p^{\rho}_{\lambda_k}(X)$ from $\mathbf{S}_{\rho}$ to $\mathbf{S}_{\lambda_k}$ is well-defined, where $X$ is in $\Gamma(M,T(M))$. Then we have a formula of $D_{\lambda_k}^{\rho}$,
%%%%%%
\begin{equation}
D^{\rho}_{\lambda_k}=\sum p^{\rho}_{\lambda_k}(e_i)\nabla_{e_i}:\Gamma(M,\mathbf{S}_{\rho})\to \Gamma(M,\mathbf{S}_{\lambda_k}).    \label{eqn:3-4}
\end{equation}
%%%%
Here, $\{e_i\}$ is an orthonormal frame of $T(M)\simeq T^{\ast}(M)$. We set a formal adjoint operator of $D^{\rho}_{\lambda_k}$ by $(D^{\rho}_{\lambda_k})^{\ast}:=\sum -(p^{\rho}_{\lambda_k}(e_i))^{\ast}\nabla_{e_i}$. 

We shall extend the relations among the Clifford homomorphisms to the ones among higher spin Dirac operators. First, we show that the sum of the higher spin Dirac operator $(D^{\rho}_{\lambda_k})^{\ast}D^{\rho}_{\lambda_k}$ for $0\le k\le N$ is the connection Laplacian.
%%%%%%%%%%%%%%%%%%%%%%%%%%%%%%%%%%%%%%%%%%%%%%%%%%%%%%%%%
%%%%%%%%%%% proposition 3-1  %%%%%%%%%%%%%%%%%%%%%%%%%%%%
%%%%%%%%%%%%%%%%%%%%%%%%%%%%%%%%%%%%%%%%%%%%%%%%%%%%%%%%
\begin{prop}\label{prop:3-1}
Let $\nabla^{\ast}\nabla$ is the connection Laplacian on $\Gamma(M,\mathbf{S}_{\rho})$ defined by $\sum -\nabla_{e_i}\nabla_{e_i}$, and $\sum \mathbf{S}_{\lambda_k}$ be the irreducible decomposition of $\mathbf{S}_{\rho}\otimes T^{\ast}(M)$. Then 
%%%%%%%%%%%%%
\begin{equation}
 \nabla^{\ast}\nabla=\sum_k (D^{\rho}_{\lambda_k})^{\ast}D^{\rho}_{\lambda_k}. 
                                  \label{eqn:3-6}
\end{equation}
%%%%%%%%%%
\end{prop}
%%%%%%%%%%%%%%%%%%%%%%%%%%%%%%%%%%%%%%%
%%%%%%%%%%%%%%%%%%%  begin proof %%%%%%%%%%%%%
\begin{proof}
For any $x$ in $M$, we can find a local orthonormal frame $\{e_i\}$ such that $(\nabla_{e_i}e_j)_x=0$ for any $i$ and $j$, so that $(\nabla_{e_i}p^{\rho}_{\lambda_k}(e_j))_x=0$. Then we have 
%%%%%%%%%
\begin{equation}
\begin{split}
\sum_k (D^{\rho}_{\lambda_k})^{\ast}D^{\rho}_{\lambda_k}&=
     \sum_{i,j}\sum_k -(p^{\rho}_{\lambda_k}(e_j))^{\ast}
     p^{\rho}_{\lambda_k}(e_i)\nabla_{e_j}\nabla_{e_i}  \\
&=\sum_{i,j}-\delta_{ij}\nabla_{e_j}\nabla_{e_i} \\
&=\sum_i-\nabla_{e_i}\nabla_{e_i}.
\end{split}\nonumber
\end{equation}
%%%
\end{proof}
%%%%%%%%%%%%%%%%%%%% end proof %%%%%%%%%%%%%%%%%%%%
Next, we introduce a curvature transformation on $\mathbf{S}_{\rho}$. For $X$ and $Y$ in $\Gamma(M,T(M))$, the curvature $R_{\rho}(X,Y)$ is defined by 
%%%%%%%%%%
\begin{equation}
R_{\rho}(X,Y)=\nabla_X \nabla_Y-\nabla_Y \nabla_X-\nabla_{[X,Y]} \in \Gamma (M,\mathrm{End}(\mathbf{S}_{\rho})).   \label{eqn:3-7}
\end{equation}
%%%%%%%%
This curvature leads the curvature transformation $R^1_{\rho}$ in $\Gamma (M,\mathrm{End}(\mathbf{S}_{\rho}))$,
%%%%%%%
\begin{equation}
R_{\rho}^1:=\frac{1}{8}\sum_{i,j} \pi_{\rho}([e_i,e_j])R_{\rho}(e_i,e_j).  
                                 \label{eqn:3-8}
\end{equation}
%%%%%%%%
For example, $R^1_{\Delta}$ is $\kappa/8$ and $R^1_{\Lambda^1}$ is $\mathrm{Ric}$, where $\kappa$ is the scalar curvature and $\mathrm{Ric}$ is the Ricci curvature transformation. We show that the sum of $(D^{\rho}_{\lambda_k})^{\ast}D^{\rho}_{\lambda_k}$ with conformal weight is the curvature transformation. 
%%%%%%%%%%%%%%%%%%%%%%%%%%%%%%%%%%%%%%%%%%%%%%%%%%%%%%%%%%%%%%%%%
%%%%%%%%%%%%%%%%%%%%%%%  proposition 3-2 %%%%%%%%%%%%%%%%%%%%%%%
%%%%%%%%%%%%%%%%%%%%%%%%%%%%%%%%%%%%%%%%%%%%%%%%%%%%%%%%%%%%%%%%
\begin{prop}\label{prop:3-2}
Let $R^1_{\rho}$ be the curvature transformation on $\mathbf{S}_{\rho}$ as above. Then we have 
%%%%%%%
\begin{equation}
 R_{\rho}^1=\sum_k m(\lambda_k)(D^{\rho}_{\lambda_k})^{\ast}D^{\rho}_{\lambda_k}.                                 \label{eqn:3-9}
 \end{equation}
%%%%%%% 
\end{prop}
%%%%%%%%%%%%%%%%%%%%%%%%%%%%%%%%%%%%%%%%%%%%%%%%%%%%%%%%%%%%

Now, we define a family of self adjoint differential operators on $\Gamma (M,\mathbf{S}_{\rho})$ by $R_{\rho}^{q}:=\sum_{i,j} -r^q_{\rho}(e_i,e_j)\nabla_i\nabla_j$ for $q\ge 0$. If $q$ is even and $M$ is compact, then $R^q_{\rho}$ is a non-negative operator and $\ker R^q_{\rho}$ is $\ker \nabla$, that is, the space of parallel sections. 
%%%%%%%%%%%%%%%%%%%%%%%%%%%%%%%%%%%%%%%%%%%%%%%%%%%%%%%%%%%%%%%%%
%%%%%%%%%%%%%%%%%%%%%%%  theorem 3-3    %%%%%%%%%%%%%%%%%%%%%%%
%%%%%%%%%%%%%%%%%%%%%%%%%%%%%%%%%%%%%%%%%%%%%%%%%%%%%%%%%%%%%%%%
\begin{thm}\label{thm:3-3}
The sum of higher spin Dirac operators $(D^{\rho}_{\lambda_k})^{\ast}D^{\rho}_{\lambda_k}$ with the $q$-th power of conformal weight is  $R_{\rho}^q$: 
%%%%%%%%%%%
\begin{equation}
 R_{\rho}^q=\sum_k m(\lambda_k)^q(D^{\rho}_{\lambda_k})^{\ast}D^{\rho}_{\lambda_k}.                \label{eqn:3-10}
 \end{equation}
 %%%%%%%%
This implies an explicit formula of $(D_{\lambda_k}^{\rho})^{\ast}D_{\lambda_k}^{\rho}$ as follows:
%%%%%%%
 \begin{equation}
(D_{\lambda_k}^{\rho})^{\ast}D_{\lambda_k}^{\rho}=\sum_{q} n_{k q}R^q_{\rho},
                          \label{eqn:3-11}
\end{equation}
%%%%%%%%%
where $n_{ij}$ is a constant given in \eqref{eqn:2-20} and, for the exceptional case, we think of the direct sum of $\mathbf{S}_{\rho+\mu_m}$ and $\mathbf{S}_{\rho-\mu_m}$ as only one component of $\sum \mathbf{S}_{\lambda_k}$. 
\end{thm}
%%%%%%%%%%%%%%%%%%%%%%%%%%%%%%%%%%%%%%
%%%%%%%%%%%%%%%%%%%%%%%%%% 
%%%%% remark 3-3 %%%%%%%%%
%%%%%%%%%%%%%%%%%%%%%%%%%%
\begin{rem}
We can define a differential operator $\mathrm{P}_{\rho}$ by $\mathrm{P}_{\rho}:=\sum \mathrm{pf}_{\rho}(e_i,e_j)\nabla_{e_i}\nabla_{e_j}$. Since $\mathrm{pf}_{\rho}(e_i,e_j)=-\mathrm{pf}_{\rho}(e_j,e_i)$ for the exceptional case, the operator $P_{\rho}$ is a curvature transformation.
\end{rem}
%%%%%%%%%%%%%%%%%%%
%%%%%%%%%%%%%%%%%%%%%%%%%%%%%%%%%%%%%%%%%%%%%%%%%%%%%%
%%%%%%%%%%% example %%%%%%%%%%%%%%%%%%%%%%%%%%%%%%%%%%
%%%%%%%%%%%%%%%%%%%%%%%%%%%%%%%%%%%%%%%%%%%%%%%%%%%%%%
\begin{ex}[Spinor]\label{ex:spinor 2}
We have two higher spin Dirac operators on spinor bundle $\mathbf{S}_{\Delta}$, that is, the Dirac operator $D$ and the twistor operator $T$. Our higher spin Dirac operators $D_{\Delta}^{\Delta}:\Gamma (M,\mathbf{S}_{\Delta})\to \Gamma (M,\mathbf{S}_{\Delta})$ and $D_{T}^{\Delta}:\Gamma (M,\mathbf{S}_{\Delta})\to \Gamma (M,\mathbf{S}_{T})$ satisfy that
%%%%%%%%%%
\begin{gather}
(D^{\Delta}_{\Delta})^{\ast}D^{\Delta}_{\Delta}+(D^{\Delta}_{T})^{\ast}D^{\Delta}_{T}=\nabla^{\ast}\nabla,  \label{eqn:3-12}  \\
-\frac{1}{2}(D^{\Delta}_{T})^{\ast}D^{\Delta}_{T}+\frac{n-1}{2}(D^{\Delta}_{\Delta})^{\ast}D^{\Delta}_{\Delta}=R_{\Delta}^1=\frac{1}{8}\kappa. \label{eqn:3-13}
\end{gather}
%%%%%%%%%
So we set $D^{\Delta}_{\Delta}:=\frac{1}{\sqrt{n}}D$ and $D^{\Delta}_{T}:=\frac{\sqrt{n-1}}{\sqrt{n}}T$. Then we have Bochner identities, 
%%%%%%%%
\begin{equation}
D^2=\nabla^{\ast}\nabla+\frac{1}{4}\kappa, \quad T^{\ast}T=\nabla^{\ast}\nabla-\frac{1}{4(n-1)}\kappa.          \label{eqn:3-14}
\end{equation}
%%%%%%
It is from the projection formula \eqref{eqn:2-e-6} that twistor operator is realized as follows (see \cite{BFGK} and \cite{BF}):
%%%%%%%%%%%%%%
\begin{equation}
D^{\Delta}_{T}(\phi)=\sum_i p^{\Delta}_{T}(e_i)\nabla_{e_i}\phi=\sum (\nabla_{e_i}\phi +\frac{1}{n}e_i D \phi)\otimes e_i. 
                                \label{eqn:3-15}
\end{equation}
%%%%%%%%%%%%%%%%%%

Now, we assume that $M$ is compact and $\lambda^0$ is the first eigenvalue of $D^2$. To give a lower bound of $\lambda^0$, we rewrite the identity \eqref{eqn:3-13} as
%%%%%%%%
\begin{equation}
D^2=\frac{n}{4(n-1)}\kappa+\frac{n-1}{2n}T^{\ast}T. \label{eqn:3-16}
\end{equation}
%%%%%%%
This identity indicates that $\lambda^0\ge \frac{n}{4(n-1)}\kappa_0$, where $\kappa_0$ is $\min_{x\in M}\kappa(x)$ (see \cite{BFGK} and \cite{BF}). 
\end{ex}
%%%%%%%%%%%%%%%%%%%%%%%%%%%%%%%%%%%%%%%%%%%%%%%%%%%%%%
%%%%%%%%%%% example %%%%%%%%%%%%%%%%%%%%%%%%%%%%%%%%%%
%%%%%%%%%%%%%%%%%%%%%%%%%%%%%%%%%%%%%%%%%%%%%%%%%%%%%%
\begin{ex}[Differential from]\label{ex:diff 2}
We consider the $k$-th differential forms $\Lambda^k(M)$. Then we obtain three higher spin Dirac operators,
%%%%%%%%
\begin{gather} 
d:=\sqrt{k+1}D_{\lambda_1}^{\rho}:\Gamma (M,\Lambda^k(M)) \to \Gamma(M,\Lambda^{k+1}(M)),  \label{eqn:3-17}\\
d^{\ast}:=-\sqrt{n-k+1}D_{\lambda_2}^{\rho}:\Gamma(M,\Lambda^k(M))\to \Gamma(M,\Lambda^{k-1}(M)),  \label{eqn:3-18}\\ 
C:=D_{\lambda_0}^{\rho}:\Gamma(M,\Lambda^k(M))\to \Gamma(M,\mathbf{S}_{\lambda_0}).           \label{eqn:3-19}
\end{gather}
%%%%%%%%%
Here $C$ is called the conformal killing operator whose kernel is the space of conformal killing forms. It follows from Theorem \ref{thm:3-3} that 
%%%%%%
\begin{gather}
C^{\ast}C+\frac{1}{k+1}d^{\ast}d+\frac{1}{n-k+1}dd^{\ast}=\nabla^{\ast}\nabla,                      \label{eqn:3-20}\\
-C^{\ast}C+\frac{k}{k+1}d^{\ast}d+\frac{n-k}{n-k+1}dd^{\ast}=R^1_{\rho}, 
                   \label{eqn:3-21}\\
C^{\ast}C+\frac{k^2}{k+1}d^{\ast}d+\frac{(n-k)^2}{n-k+1}dd^{\ast}=R^2_{\rho}.
                     \label{eqn:3-22}
\end{gather}
%%%%
The first and second identities give the Bochner identity for differential forms. The projection formula \eqref{eqn:2-e-15} leads a formula of $C$ (see \cite{K}): 
%%%%%%%
\begin{equation}
C(\phi)=\sum (\nabla_i\phi-\frac{1}{k+1}i(e_i) d\phi+ \frac{1}{n-k+1}e_{i\wedge}d^{\ast}\phi) \otimes e_i.  \label{eqn:3-23}
\end{equation}
%%%%%%
\end{ex}
%%%%%%%%%%%%%%%%%%%%%%%%%%%%%%%%%%%%%%%%%%%%%%%%%%%%%%
%%%%%%%%%%% example %%%%%%%%%%%%%%%%%%%%%%%%%%%%%%%%%%
%%%%%%%%%%%%%%%%%%%%%%%%%%%%%%%%%%%%%%%%%%%%%%%%%%%%%%
\begin{ex}[Higher spin fields on anti self-dual $4$-manifolds] 
Let $M$ be an anti self-dual $4$-dimensional spin manifold and $\mathbf{S}_{k,l}$ be the vector bundle corresponding to the highest weight $\rho:=(k+l,k-l)$ for $k,l$ in $\mathbf{Z}\cup(\mathbf{Z}+1/2)$. We remark that $\mathbf{S}_{k,l}$ corresponds to the representation of $\pi_k\widehat{\otimes}\pi_l$ of $Spin(4)=SU(2)\times SU(2)$, where $\pi_k$ is a spin $k$ representation of $SU(2)$. We consider the vector bundle $\mathbf{S}_{k,0}$ and have two higher spin Dirac operators,
%%%
\begin{gather}
D_0:\Gamma (M,\mathbf{S}_{k,0})\to \Gamma (M,\mathbf{S}_{k+1/2,1/2}),
               \label{eqn:3-24}\\
D_1:\Gamma (M,\mathbf{S}_{k,0})\to \Gamma (M,\mathbf{S}_{k-1/2,1/2}). 
               \label{eqn:3-25}
\end{gather}
%%%%%
The relations between $D_0$ and $D_1$ are
%%%%%%
\begin{gather}
(D_0)^{\ast}D_0+(D_1)^{\ast}D_1=\nabla^{\ast}\nabla, \label{eqn:3-26}\\
-k(D_0)^{\ast}D_0+(k+1)(D_1)^{\ast}D_1=R_{k,0}^1.  \label{eqn:3-27}
\end{gather}
%%%%%%%
Since $M$ is an anti self-dual manifold, we can show that the curvature transformation $R_{k,0}^1$ is $\frac{k(k+1)}{6}\kappa$. We set $\widehat{D}_1:=\sqrt{(2k+1)/k}D_1$ and give the Bochner identity for $\widehat{D}_1$ (see \cite{AHS} and \cite{H}):
%%%%%%%
\begin{equation}
(\widehat{D}_1)^{\ast}\widehat{D}_1=\nabla\nabla^{\ast}+\frac{k+1}{6}\kappa.
                    \label{eqn:3-28}
\end{equation}
%%%%%%%
The identity \eqref{eqn:3-27} gives a lower bound of the first eigenvalue $\lambda^0$ of $(\widehat{D}_1)^{\ast}\widehat{D}_1$, that is, $\lambda^0\ge \frac{2k+1}{6}\kappa_0$.
\end{ex}
%%%%%%%%%%%%%%%%%%%%%%%%%%%%%%%%
%%            4               %%
%%%%%%%%%%%%%%%%%%%%%%%%%%%%%%%%
\section{General Bochner identities} \label{sec:4}
We have considered some Bochner identities in the previous section. It allows us to give a general Bochner identity on any associated bundle $\mathbf{S}_{\rho}$. We have known the following two identities: 
%%%%%%
\begin{gather}
\sum_{0\le k\le N} (D^{\rho}_{\lambda_k})^{\ast}D^{\rho}_{\lambda_k}=\nabla^{\ast}\nabla ,                    \label{eqn:4-1}  \\
\sum_{0\le k\le N} m(\lambda_k)(D^{\rho}_{\lambda_k})^{\ast}D^{\rho}_{\lambda_k}=R_{\rho}^1,                  \label{eqn:4-2}
\end{gather}
%%%%%
where we order $\{\lambda_k\}$ to satisfy that $\lambda_0>\lambda_1>\cdots>\lambda_N$. We eliminate $(D^{\rho}_{\lambda_0})^{\ast}D^{\rho}_{\lambda_0}$ from above two equations and obtain a general Bochner identity on $\mathbf{S}_{\rho}$,
%%%
\begin{equation}
\sum_{1\le k\le N}(1-\frac{m(\lambda_k)}{m(\lambda_0)})(D^{\rho}_{\lambda_k})^{\ast}D^{\rho}_{\lambda_k}=\nabla^{\ast}\nabla+\frac{1}{-m(\lambda_0)}R_{\rho}^1.
                              \label{eqn:4-3}
\end{equation}
%%%%
Here, $m(\lambda^0)$ is $-\rho^1$ and $(1-\frac{m(\lambda_k)}{m(\lambda_0)})$ is positive. 
%%%%%%%%%%%%%%%%%%%%%%%%%%%%%%%%%%%%%%%%%%%%%%%%%%%%%%%%%%
%%%%%%%%%%%%%%   theorem  4-1 %%%%%%%%%%%%%%%%%%%%%%%%%%%%
%%%%%%%%%%%%%%%%%%%%%%%%%%%%%%%%%%%%%%%%%%%%%%%%%%%%%%%%%%
\begin{thm}\label{thm:4-1}
We define the differential operator $\widehat{D}^{\rho}_{\lambda_k}$ for $1\le k\le N$ by 
%%%
\begin{equation}
\widehat{D}^{\rho}_{\lambda_k}:=\sqrt{1-\frac{m(\lambda_k)}{m(\lambda_0)}}D^{\rho}_{\lambda_k}.                 \label{eqn:4-4}
\end{equation}
%%%
Then the following identity holds:
%%%%%
\begin{equation}
\Delta_{\rho}:=\sum_{k\ge 1}(\widehat{D}^{\rho}_{\lambda_k})^{\ast}\widehat{D}^{\rho}_{\lambda_k}=\nabla^{\ast}\nabla+\frac{1}{\rho^1}R_{\rho}^1.  
                               \label{eqn:4-5}
\end{equation}
%%%%%%
This Laplace type operator $\Delta_{\rho}$ is an elliptic second order operator on $\Gamma (M,\mathbf{S}_{\rho})$. If $M$ is compact, then $\Delta_{\rho}$ is non-negative and $\ker \Delta_{\rho}$ is the intersection of $\ker \widehat{D}^{\rho}_{\lambda_k}$ for $1\le k\le N$.  
\end{thm}
%%%%%%%%%%%%%%%%%%%%%%%%%%%%%%%%%%%%%%%%%%%%%%%%%%
%%%%%%%%%%%%%%%%%%%%%%%%%%%%%%%%
%%            5               %%
%%%%%%%%%%%%%%%%%%%%%%%%%%%%%%%%
\section{A lower bound of $\Delta_{\rho}$ on manifolds of positive curvature}\label{sec:5}
We assume that $M$ is a compact manifold of positive curvature in this section. In other words, there exists a constant $r>0$ such that the curvature $R_{\Ad}(\cdot,\cdot)$ on the tangent bundle $T(M)$ satisfies 
%%%%%%%%%%
\begin{equation}
\langle R_{\Ad}(u,v)v,u \rangle \ge 2r\langle u,u\rangle \langle v,v \rangle \quad 
      \textrm{for any $u$ and $v$ in $T(M)$ }. \label{eqn:5-1}
 \end{equation}
 %%%%%%%%%%%
Then the curvature transformation $R_{\rho}^1$ on $\mathbf{S}_{\rho}$ has a lower bound as follows: for $\phi$ in $\mathbf{S}_{\rho}$,  
%%%%%%%%%%%
\begin{equation}
\begin{split}
 \langle R_{\rho}^1 \phi,\phi\rangle &=\frac{1}{4} \langle \sum_{i,j} \pi_{\rho}([e_i,e_j])R_{\rho}(e_i,e_j)\phi,\phi \rangle  \\
   &=\frac{1}{64}\sum_{i,j,k,l} \langle R_{\Ad}(e_i,e_j)e_k,e_l \rangle \langle \pi_{\rho}([e_i,e_j])\pi_{\rho}([e_k,e_l])\phi, \phi \rangle \\ 
  &\ge \frac{r}{32}\sum \langle  \pi_{\rho}([e_i,e_j])\pi_{\rho}([e_j,e_i])\phi,\phi \rangle \\
 &=-2rc(\rho)\langle \phi,\phi \rangle, 
\end{split}    \label{eqn:5-2}
\end{equation}
%%%%%%%%%
where $c(\rho)$ is a negative constant due to the Casimir operator on $V_{\rho}$ in \eqref{eqn:1-9}.
%%%%%%%%%%%%%%%%%%%%%%%%%%%%%%%%%%%%%%%%%%%%%%%
%%%%%%%%%%%%%%    theorem   %%%%%%%%%%%%%%%%%%%
%%%%%%%%%%%%%%%%%%%%%%%%%%%%%%%%%%%%%%%%%%%%%%%
\begin{thm}\label{thm:5-1}
Let $M$ be a compact manifold of positive curvature as above and $\Delta_{\rho}$ be the Laplace type operator on $\Gamma (M,\mathbf{S}_{\rho})$ given in \eqref{eqn:4-5}. The first eigenvalue $\lambda^0$ of $\Delta_{\rho}$  has a lower bound, 
%%%
\begin{equation}
\lambda^0 \ge \frac{m(\lambda_0)-m(\lambda_N)}{m(\lambda_0)m(\lambda_N)}(-2rc(\rho))>0.  \label{eqn:eigen}
\end{equation}
%%%%%
If the equality holds in \eqref{eqn:eigen}, then the eigensections with the first eigenvalue $\lambda^0$ are in $\ker D_{\lambda_0}^{\rho}$ and $\ker \widehat{D}_{\lambda_k}^{\rho}$ for $1\le k \le N-1$. 
\end{thm}
%%%%%%%%%%%%%%%%%%%%%
%%%%%%%%%%%%%%%%  begin proof %%%%%%%%%%%%%%%%%%%%%%%%%
\begin{proof}
We remark that $m(\lambda_N)>0$ and, for any $k$,
%%%%%%
\begin{equation}
\frac{m(\lambda_0)m(\lambda_N)}{m(\lambda_0)-m(\lambda_N)}\ge \frac{m(\lambda_0)m(\lambda_k)}{m(\lambda_0)-m(\lambda_k)}. \nonumber
\end{equation}
%%%
We have known the following identities: 
%%%
\begin{gather}
\Delta_{\rho}=\sum (\widehat{D}^{\rho}_{\lambda_k})^{\ast}\widehat{D}^{\rho}_{\lambda_k},  \nonumber \\
 R_{\rho}^1+\rho^1 (D_{\lambda_0}^{\rho})^{\ast}D_{\lambda_0}^{\rho}=\sum_{1\le k\le N} \frac{m(\lambda_0)m(\lambda_k)}{m(\lambda_0)-m(\lambda_k)}(\widehat{D}^{\rho}_{\lambda_k})^{\ast}\widehat{D}^{\rho}_{\lambda_k}. \nonumber
\end{gather}
%%%
Then 
%%%%%%
\begin{equation}
\begin{split}
(\Delta_{\rho} \phi,\phi)&=\sum (\widehat{D}^{\rho}_{\lambda_k}\phi,\widehat{D}^{\rho}_{\lambda_k}\phi) \\
&=\frac{m(\lambda_0)-m(\lambda_N)}{m(\lambda_0)m(\lambda_N)}\sum \frac{m(\lambda_0)m(\lambda_N)}{m(\lambda_0)-m(\lambda_N)} \|\widehat{D}^{\rho}_{\lambda_k}\phi\|^2 \\
&\ge \frac{m(\lambda_0)-m(\lambda_N)}{m(\lambda_0)m(\lambda_N)}\sum \frac{m(\lambda_0)m(\lambda_k)}{m(\lambda_0)-m(\lambda_k)} \|\widehat{D}^{\rho}_{\lambda_k}\phi\|^2 \\
&=\frac{m(\lambda_0)-m(\lambda_N)}{m(\lambda_0)m(\lambda_N)}((R_{\rho}^1\phi,\phi)+\rho^1 \|D_{\lambda_0}^{\rho}\phi \|^2 )\\
&\ge \frac{m(\lambda_0)-m(\lambda_N)}{m(\lambda_0)m(\lambda_t)}(R_{\rho}^1\phi,\phi)=\frac{m(\lambda_0)-m(\lambda_N)}{m(\lambda_0)m(\lambda_N)}\int_M \langle R_{\rho}^1\phi,\phi\rangle \ast 1\\
&\ge \frac{m(\lambda_0)-m(\lambda_N)}{m(\lambda_0)m(\lambda_N)}(-2rc(\rho))(\phi,\phi).
\end{split}\nonumber
\end{equation}
This inequality leads us to the proposition.
\end{proof}
%%%%%%%%%%%%%%%%% end proof %%%%%%%%%%%%%%%
%%%%%%%%%%%%%%%%%%%%%%%%%%%%%%%%%%%%%%%%%%%%%%%%%
%%%%%%%%%% example %%%%%%%%%%%%%%%%%%%%%%%%%%%%%
%%%%%%%%%%%%%%%%%%%%%%%%%%%%%%%%%%%%%%%%%%%%%%%%%
\begin{ex}[Differential form]\label{ex:diff 3}
We give a lower bound of the first eigenvalue $\lambda^0$ of $dd^{\ast}+d^{\ast}d$ on $\Lambda^k(M)$, where $M$ is a compact $n$-dimensional manifold of positive curvature. From the above theorem, we have $\lambda^0\ge k(n-k+1)r$ for $1\le k \le [\frac{n}{2}]$. If the equality holds in the equation, the eigensections with the eigenvalue $\lambda^0$ is in $\ker C\cap \ker d$. This lower bound coincides with the one given in \cite {GM}.
\end{ex}
%%%%
%%%%%%%%%%%%%%%%%%%%%%%%%%%%%%%%%%%%%%%
%%  acknowledgements            %%%%%%%
%%%%%%%%%%%%%%%%%%%%%%%%%%%%%%%%%%%%%%%
\section*{Acknowledgements}
The author is partially supported by Waseda University Grant for Special Research Project 2000A-880.
%%%%%%%%%%%%%%%%%%%%%%%%%%%%%%%%
%%           6                %%
%%%%%%%%%%%%%%%%%%%%%%%%%%%%%%%%
\appendix
\section{Appendix: Conformal invariance of the higher spin Dirac operator}\label{sec:6}
In \cite{F}, Fegan shows that all the higher spin Dirac operators are conformal invariant first order differential operators. We show the conformal invariance of $D^{\rho}_{\lambda}$ explicitly by using the Clifford homomorphisms. Let $(M,g)$ be a spin manifold with Riemannian metric $g$. We deform the metric $g$ conformally to $g':=e^{2\sigma(x)} g$, where $\sigma(x)$ is a scalar function on $M$. We denote the objects associated to $g'$ by adding a symbol `` $'$ " to them: for example, $\mathbf{Spin}'(M)$, $\mathbf{S}'_{\rho}$ and so on. 

The isomorphism between the orthonormal frame bundles $\mathbf{SO}(M)$ and $\mathbf{SO}'(M)$ is realized by the mapping 
%%%%%%%%%%
\begin{equation}
\Psi: \mathbf{SO}(M)\ni \{e_i\}_i\mapsto \{e_i':=e^{-\sigma}e_i\}_i\in \mathbf{SO}'(M).    \label{eqn:6-1}
\end{equation}
%%%%%%%%
This mapping is lifted to the isomorphism for spin bundles and induces a bundle isometry $\psi_{\rho}:\mathbf{S}_{\rho}\to\mathbf{S}'_{\rho}$ for each associated bundle such that $\psi_{\lambda}\circ p^{\rho}_{\lambda}(e_i)=e^{\sigma}p^{\rho}_{\lambda}(e_i)\circ \psi_{\rho}$. The covariant derivatives $\nabla$ and $\nabla'$ on $\Gamma(M,\mathbf{S}_{\rho})$ and $\Gamma(M, \mathbf{S}'_{\rho})$ are related as follows (for the spinor case, see \cite{BFGK} and \cite{LM}):
%%%%%%
\begin{equation}
\nabla'_X=\psi_{\rho}\circ \{\nabla_X\phi +\frac{1}{4}\pi_{\rho}([\mathrm{grad} (\sigma),X])\}\circ\psi_{\rho}^{-1},   \label{eqn:6-2}
\end{equation}
%%%%%%%
where $X$ is any vector field on $M$. This relation implies a relation for the higher spin Dirac operator.
%%%%%%%%%%%%%%%%%%%%%%%%%%%%%%%%%%%%%%%%%%%%%%%%%%%%%%%%%%%%%
%%%%%%%%%%  lemma %%%%%%%%%%%%%%%%%%%%%%%%%%%%%%%%%%%%%%%%%
%%%%%%%%%%%%%%%%%%%%%%%%%%%%%%%%%%%%%%%%%%%%%%%%%%%%%%%%%%
\begin{lem}\label{lem:6-1}
The higher spin Dirac operators $D^{\rho}_{\lambda}$ and $D'{}^{\rho}_{\lambda}$ associated to metric $g$ and $g'$ respectively are related as follows:
%%%%%%
\begin{equation}
D'{}^{\rho}_{\lambda}=e^{-\sigma}\psi_{\lambda} \circ\{D^{\rho}_{\lambda}+m(\lambda)p^{\rho}_{\lambda}(\mathrm{grad} (\sigma))\}\circ\psi_{\rho}^{-1}         \label{eqn:6-3}
\end{equation}
\end{lem}
%%%%%%%%%%%%%%%%%%%%%%%%%%%%%%
%%%%%%%%%%%%%  begin proof %%%%%%%%%%%%%%%%%%%%%%%%%
\begin{proof}
For $\phi$ in $\Gamma (M,\mathbf{S}_{\rho})$, we have
\begin{equation}
\begin{split}
D'{}^{\rho}_{\lambda}\circ\psi_{\rho}(\phi)&=\sum p^{\rho}_{\lambda}(e'_i)\nabla'_{e'_i}\psi_{\rho}(\phi) \\
&=\sum p^{\rho}_{\lambda}(e_i)\nabla'_{e_i}\circ\psi(\phi) \\
&=\sum p^{\rho}_{\lambda}(e_i) \psi_{\rho} \circ \{ \nabla_{e_i}\phi +\frac{1}{4}\pi_{\rho}([\mathrm{grad} (\sigma),e_i])\phi\} \\
&=e^{-\sigma}\psi_{\lambda}\circ\{ \sum p^{\rho}_{\lambda}(e_i)(\nabla_{e_i}\phi +\frac{1}{4}\pi_{\rho}([\mathrm{grad} (\sigma),e_i])\phi) \}  \\
&=e^{-\sigma}\psi_{\lambda}\circ\{ D^{\rho}_{\lambda} \phi + m(\lambda)p^{\rho}_{\lambda}(\mathrm{grad} (\sigma))\phi  \}  \quad \textrm{( from \eqref{eqn:2-7})} .
\end{split}\nonumber
\end{equation}
\end{proof}
%%%%%%%%%%%%% proof end %%%%%%%%%
The usual Dirac operator $D$ satisfies that $D\circ f =\mathrm{grad}(f)\cdot+fD$, where $f$ is any smooth function on $M$. The higher spin Dirac operator has the similar property. 
%%%%%%%%%%%%%%%%%%%%% proof end %%%%%%%%
%%%%%%%%%%%%%%%%%%%%%%%%%%%%%%%%%%%%%%%%%%%%%%%%%%%%%%%%%%%%%
%%%%%%%%%%  lemma %%%%%%%%%%%%%%%%%%%%%%%%%%%%%%%%%%%%%%%%%
%%%%%%%%%%%%%%%%%%%%%%%%%%%%%%%%%%%%%%%%%%%%%%%%%%%%%%%%%%
\begin{lem}\label{lem:6-2}
For any smooth function $f$ on $M$, we have 
%%%%%%%
\begin{equation}
D^{\rho}_{\lambda}\circ f =p^{\rho}_{\lambda}(\mathrm{grad}(f))+fD^{\rho}_{\lambda}.                          \label{eqn:6-4}
\end{equation}
%%%%%%
\end{lem}
%%%%%%%%%%%%%%%%%%%%%%%%%%%%
The above two lemma gives us the conformal invariance of the operator $D^{\rho}_{\lambda}$:
%%%%%%%%%%%%%%
\begin{equation}
D'{}^{\rho}_{\lambda}=(e^{-(m(\lambda)+1)\sigma}\psi_{\lambda})\circ D^{\rho}_{\lambda}\circ (e^{-m(\lambda)\sigma}\psi_{\rho})^{-1}.     \label{eqn:6-5}
\end{equation}
%%%%%%%%%%%%%%%%%%
In particular, we have $\dim \ker D'{}^{\rho}_{\lambda}=\dim \ker D^{\rho}_{\lambda}$.
%%%%%%%%%%%%%%%%%%%%%%%%%%%%%%%%
%        reference             %
%%%%%%%%%%%%%%%%%%%%%%%%%%%%%%%%
%%

%%%%%%%%%%%%%%%%%%%%%%%%%
\end{document}